\documentclass[a4paper,12pt,leqno]{amsart}
\usepackage{color}
\usepackage{amsfonts}
\usepackage{amssymb}
\usepackage{amsmath}
\usepackage{amsthm}
\usepackage{hyperref}
\usepackage{mathtools}
\usepackage{enumitem}
\usepackage{subcaption}
\usepackage{dsfont}
\usepackage{scalerel,stackengine}
\stackMath
\newcommand\reallywidehat[1]{%
\savestack{\tmpbox}{\stretchto{%
  \scaleto{%
    \scalerel*[\widthof{\ensuremath{#1}}]{\kern-.6pt\bigwedge\kern-.6pt}%
    {\rule[-\textheight/2]{1ex}{\textheight}}
  }{\textheight}%
}{0.5ex}}%
\stackon[1pt]{#1}{\tmpbox}%
}
\usepackage{bbm}

\newcommand{\Int}{\mathop{}\!\mathrm{Int}}
\makeatletter
\DeclareRobustCommand\widecheck[1]{{\mathpalette\@widecheck{#1}}}
\def\@widecheck#1#2{%
    \setbox\z@\hbox{\m@th$#1#2$}%
    \setbox\tw@\hbox{\m@th$#1%
       \widehat{%
          \vrule\@width\z@\@height\ht\z@
          \vrule\@height\z@\@width\wd\z@}$}%
    \dp\tw@-\ht\z@
    \@tempdima\ht\z@ \advance\@tempdima2\ht\tw@ \divide\@tempdima\thr@@
    \setbox\tw@\hbox{%
       \raise\@tempdima\hbox{\scalebox{1}[-1]{\lower\@tempdima\box
\tw@}}}%
    {\ooalign{\box\tw@ \cr \box\z@}}}
\makeatother

\newcommand{\Comment}[1]{}

\def\R{\mathbb{R}}

\newcommand{\bv}{\mathbf{v}}

\newtheorem{theorem}{Theorem}[section]

\newtheorem{prop}[theorem]{Proposition}
\theoremstyle{definition}
\newtheorem{remark}[theorem]{Remark}

\numberwithin{equation}{section}

\DeclareMathOperator{\rank}{rank}

\def\R{\mathbb R}

\def\bv{\big|}
\def\beas{\begin{eqnarray*}}
\def\eeas{\end{eqnarray*}}
\def\sph{\mathbb S^{d-1}}
\def\rank{\hbox{rank}}
\def\corank{\hbox{corank}}

\newcounter{vremennyj}

\begin{document}

\author{Jos\'e Gaitan}
\address[J.~ Gaitan]{Department of Mathematics, Virginia Tech, Blacksburg, VA 24061}
\email{jogaitan@vt.edu}
\author{Allan Greenleaf}
\address[A.~ Greenleaf]{Department of Mathematics, Univerity of Rochester, Rochester, NY 14627}
\email{{allan@math.rochester.edu}}
\author{Eyvindur Ari Palsson}
\address[E.~ Palsson]{Department of Mathematics, Virginia Tech, Blacksburg, VA 24061}
\email{palsson@vt.edu}
\author{Georgios Psaromiligkos}
\address[G.~ Psaromiligkos]{Department of Mathematics, Virginia Tech, Blacksburg, VA 24061}
\email{psaromil@vt.edu}

\title[]{On restricted Falconer distance sets}
\maketitle

\begin{abstract}
We introduce a  class of Falconer distance problems, which we call of restricted type, lying between the classical version and its pinned variant. 
Prototypical restricted distance sets are the diagonal distance sets,
$k$-point configuration sets given by
$$\Delta^{diag}(E)= \{ \,|(x,x,\dots,x)-(y_1,y_2,\dots,y_{k-1})| : x, y_1, \dots,y_{k-1} \in E\, \}$$
for a compact $E\subset\R^d$ and $k\ge 3$.
We show that $\Delta^{diag}(E)$
has non-empty interior if the Hausdorff dimension of $E$ satisfies
\begin{equation}\label{eqn k condition in abstract}
\dim(E) > 
\begin{cases}
\frac{2d+1}3, & k=3 \\
\frac{(k-1)d}k,& k\ge 4.
\end{cases}
\end{equation}

We prove an extension of this to $C^\omega$ Riemannian metrics $g$  close to the
product of Euclidean metrics.
For product metrics this follows from known results on pinned distance sets, 
 but to obtain a result for general perturbations $g$ we present a sequence of proofs of partial results,
 leading up to the proof of the full result, which is
  based on estimates for multilinear Fourier integral operators.
\end{abstract}

\section{The Falconer distance problem and its many variants}

The Falconer distance problem, a continuous analogue of the celebrated Erdős distance problem asks: How large does $\dim(E)$, for a compact set $E\subseteq \mathbb{R}^d$, need to be to ensure that the Lebesgue measure of its distance set
$$ \Delta(E)= \{|x-y|,\, x,y \in E\}$$
 is positive? Here and below $\dim(E)$ denotes the Hausdorff dimension of the set $E$.  Falconer introduced this problem in 1985 in \cite{Falconer} and established the dimensional threshold $\dim (E) > \frac{d+1}{2}$. 

Further, Falconer conjectured the threshold is $\dim(E) > \frac{d}{2}$ and showed the result could not hold true strictly below that threshold.
Falconer's problem has  stimulated much activity and been the focus of many outstanding results, 
e.g., \cite{Bourgain,Wolff1,Erdogan,DuGuthOuWangWilsonZhang,DuZhang,GuthIosevichOuWang,DuIosevichOuWangZhang}.

For two compact sets $E,F \subseteq \mathbb{R}^d$ one can also consider an asymmetric version, given by
$$\Delta(E,F) := \{ |x-y|:\, x\in E, y\in F \},$$
so that $\Delta(E,E)=\Delta(E)$. 
Note all the standard proofs adapt to this setting and the threshold condition can be replaced by a lower bound on $(\dim(E) +\dim(F))/2$.

Yet another variant of the Falconer  problem was introduced by Mattila and Sj\"olin \cite{MattilaSjolin}, who asked how large does $\dim(E)$ need to be in order to ensure that $\Delta(E)$ satisfies the stronger condition of having nonempty interior, 
and showed that $\dim(E) > \frac{d+1}{2}$ is sufficient.

Both the Falconer  and Mattila-Sj\"olin problems have {\it pinned} versions,  asking how large does $\dim(E)$
need to be to guarantee that there exists an $x$ such that the pinned distance set, 
$$ \Delta^x(E):= \{|x-y|\, :\,  y \in E\},$$ 
has positive Lebesgue measure or nonempty interior. 
Peres and Schlag \cite{PeresSchlag} showed that this holds for  $\dim(E) > \frac{d+2}{2},\, d\ge 3$; 
see \cite{IosevichLiu1} for some improvements and generalization.
More recently,  improvements to thresholds in the Falconer's distance problem automatically transfer over to the pinned setting 
due to the magical formula of Liu \cite{Liu}.

Nowadays one can view the original result of Falconer as well as this one of Mattila and Sj\"olin through the same lens; 
see, e.g., \cite{Mattila2}.  As with Falconer's original problem, this has
led to considerable further work in more general settings 
\cite{IosevichMourgoglouTaylor,GreenleafIosevichTaylor1,GreenleafIosevichTaylor2,RomeroAcostaPalsson1,KohPhamandShen,
RomeroAcostaPalsson2,GreenleafIosevichTaylor3}.

\section{A new problem and motivation}

In this paper we introduce new variants of the Falconer  and Mattila-Sj\"olin problems, which we call {\it restricted}  distance problems.
\footnote{After the original version of this preprint was posted, Borges, Iosevich and Ou posted  \cite{BorIosOu1},
which also discusses restricted distance problems and in some cases obtains lower thresholds than we obtain here.
See \S\S\ref{subsec relation} for a discussion.
However, we believe that Thm. \ref{main three} is not currently accessible to the methods of  \cite{BorIosOu1}, 
and in any case the techniques used to proved it indicate 
that positive results for restricted Mattila-Sj\"olin type problems can be proven in great generality.}
These lie between the original distance problems and and their pinned variants, and when stated in general encapsulate both of them.

For a compact set $E \subseteq \mathbb{R}^d$, let $F\subseteq \mathbb{R}^{d}$ be a compact set which might depend on $E$. Defining the restricted distance set, 
$$\Delta^F(E) := \{ |x-y|\, : \, x\in F, \,y \in E\},$$ 
we  ask what lower bounds on $\dim (E)$  guarantee that $\Delta^F(E)$ has positive Lebesgue measure or  
nonempty interior.
Note that if $F$ has no dependence on $E$ then   $\Delta^F(E)=\Delta(E,F)$ and so we are in the asymmetric setting of the Falconer distance problem. 
The simplest case of a set $F$ which is dependent on $E$ is when $F=E$, so  that $\Delta^F(E)=\Delta(E)$.

 If $F=\{x_0\}$ for some  point $x_0$, fixed in advance, 
then this is similar to a pinned distance problem, but 
stronger than  the usual one, since the pin is fixed. 
(A result giving nonempty interior for the set of volumes of parallelepipeds generated by an arbitrary $x_0$ and all $d$-tuples of points in $E\subset\R^d$
is in \cite[~Thm.1.2]{GreenleafIosevichTaylor2}, where this is referred to as a {\it strongly} pinned result.)

To illustrate the type  of Falconer-type distance problems in which we are interested, 
we focus on a prototype  lying between the original and pinned versions of the distance problem.
For a compact set $E\subseteq \mathbb{R}^d$,  
let $F=\{\,(x,x)\,:\, x\in E\,\}\subset \R^{2d}$, the diagonal of $E\times E$. 
With $|\, \cdot \,|$ denoting the Euclidean norm on $\R^{2d}$,
we   ask what lower bound on $\dim(E)$ ensures that 
\begin{equation}\label{eqn diag sec 2}
\Delta^F(E\times E)= \{ \,\left|\,(x,x)-(y_1,y_2)\,\right| : x, y_1, y_2\, \in E \},
\end{equation}
which we will also denote by $\Delta^{diag}(E)$,
has positive Lebesgue measure or nonempty interior.  See Fig. \ref{fig one} below.

As noted in \cite{BorIosOu1}, in order to make the problem more interesting,   
in \eqref{eqn diag sec 2} one should impose a condition $y_1\ne y_2$  because, 
if $y_1=y_2$ were allowed, then $\Delta^{diag}(E)\supset \sqrt2\cdot\Delta(E)$,
which would then have positive Lebesgue measure or nonempty interior if $\dim(E)$ is greater than the thresholds in $\R^d$  for the standard Falconer or 
Mattila-Sj\"olin  distance problems, resp. We thus include this condition and its extensions in the statements below.
So, we define 
$$\mathring{\Delta}^{diag}(E):= \{ \,\left|\,(x,x)-(y_1,y_2)\,\right| : x, y_1, y_2\, \in E,\, y_1\ne y_2\, \};$$
for  an $l$-fold Cartesian product of a $E\subset\R^d$, and an  $F\subset\R^{ld}$ define
$$
\mathring{\Delta}^F(E):=\left\{ \, |x-y|\, :\, x\in F,\, y=(y_1,\dots,y_l)\in E^l,\, y_i\ne y_j,\, \forall i\ne j\, \right\},
$$
where $|\, \cdot \, |$ is the Euclidean norm on $\R^{ld}$.

We are now ready to pose the following set of questions generalizing the prototype:
\medskip

\noindent {\bf Restricted Falconer and Mattila-Sj\"olin Problems.} Fix $l\in\mathbb N$ and a map $\mathcal F$ 
from the collection $\mathcal C\left(\R^d\right)$ 
of compact sets in $\R^d$  to $2^{\mathcal C\left(\R^{ld}\right)}$, denoting the image of a compact $E$ by $\mathcal F_E$.

\noindent{\bf Q.} What lower bounds on $\dim(E)$  ensure that either

(i) there exists an $F \in \mathcal{F}_E$ such that $\mathring{\Delta}^F(E)$ has positive Lebesgue measure (or nonempty interior); or

(ii)  for a.e. $F\in \mathcal{F}_E$ (with respect to some measure on $\mathcal F$), or

(iii) for every $F\in \mathcal{F}_E$, 

\noindent the same property holds.

\medskip

\noindent{\bf Remarks. 1.}  
For  $l=1$, case  (i) and positive Lebesgue measure, the choice of 
$\mathcal{F}_E=\{E\}$ yields the classical Falconer distance problem, while
$\mathcal{F}_E=\{ \{x\}:\,  x\in E\}$ yields its standard pinned variant. 
On the other hand, $\mathring{\Delta}^{diag}(E)$  corresponds to $l=2$, $\mathcal F_E=\{F\}$, where $F$ is the diagonal of $E\times E$ in $\R^{2d}$.
In general, if $\mathcal F$ is a singleton, the  questions (i), (ii) or (iii) collapse into one, concerning a $3$-point configuration problem of either Falconer or
Mattila-Sj\"olin type; in this paper we will focus on the latter for $\mathring{\Delta}^{diag}(E)$ and its  $k$-point configuration generalizations.

{\bf 2.} $\R^d$ can be replaced by a smooth $d$-dimensional manifold with a smooth density, and Thm. \ref{main three} below  is formulated in this setting.


\begin{figure}
    \centering
    \includegraphics[scale=0.42]{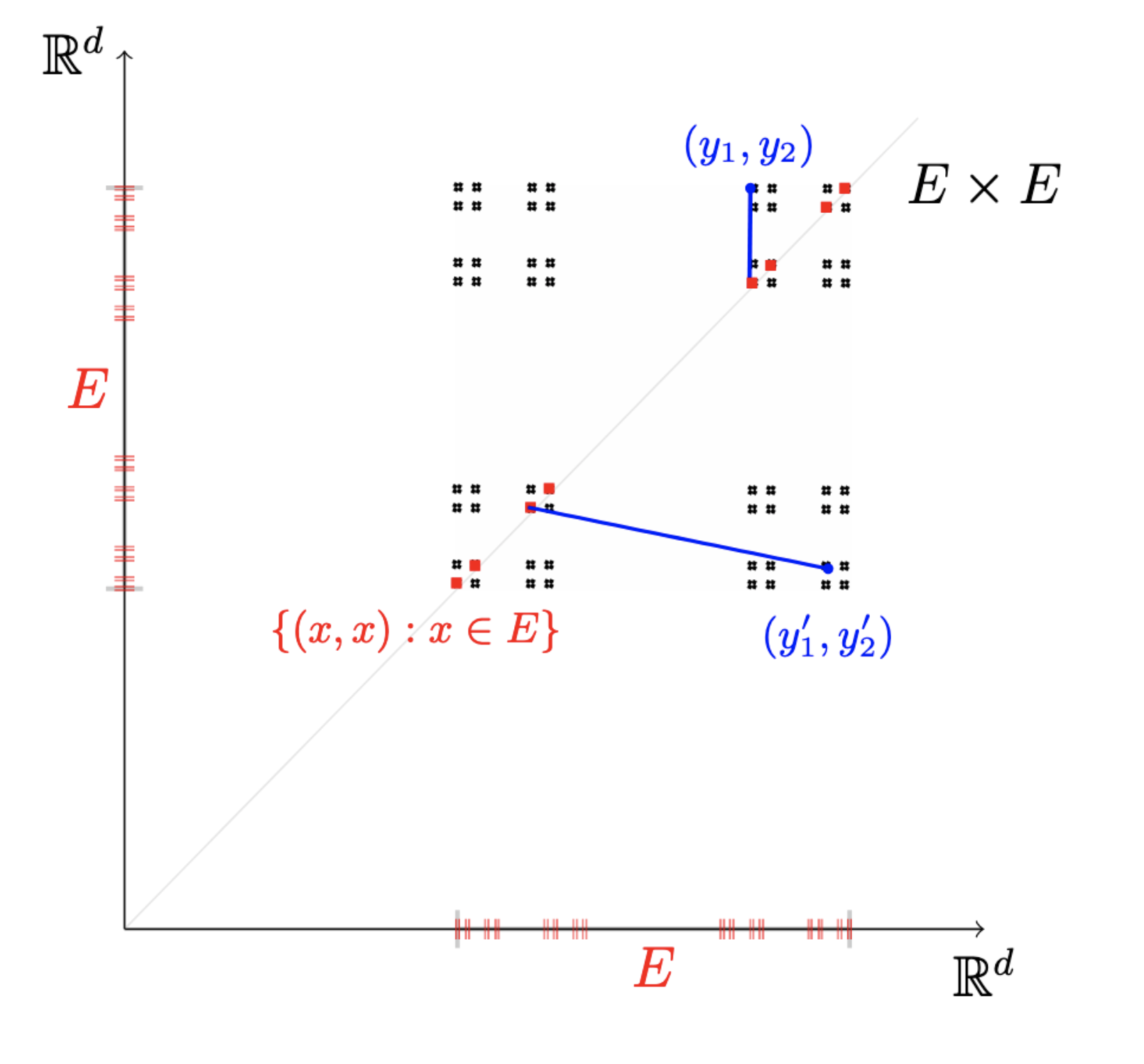}
    \caption{A sketch of how one could view $\Delta^{diag}(E)$.}
    \label{fig one}
\end{figure}


{\bf 3.} Returning to the prototype \eqref{eqn diag sec 2},
note that if we don't restrict to the diagonal but instead consider the full $\R^{2d}$ distance set $\Delta (E\times E)$, 
the best results known for the $\R^{2d}$ Falconer  problem would yield a sufficient lower bound on $\dim (E)$. 
Since $\dim(E)>\frac{d}{2}+ \frac{1}{8}$ implies that $\dim(E\times E)> \frac{2d}2+\frac14$,
and  $2d$ is even,
 the results of  \cite{GuthIosevichOuWang,DuIosevichOuWangZhang}  
yield that   $\Delta(E\times E)$ has positive Lebesgue measure.
However, for the restricted Falconer problem we are considering,  
the set $\mathring{\Delta}^{diag}(E)$ consists only of distances from points on the  diagonal of $E$ to general points of $E\times E$. 

{\bf 4.} By a result of Peres and Schlag \cite{PeresSchlag}, if $\dim(E)>(d+2)/2$, with $d\ge 3$, then there exists an $x$ such that 
the pinned distance set $\Delta^x(E)=\{\, |x-y_1|\, :\, y_1\in E\, \}$ contains an interval. 
This immediately implies that $\mathring{\Delta}^{diag}(E)$ contains an interval, since $y_2$ in \eqref{eqn diag sec 2} can simply be fixed.
The same principle applies to any $\mathring{\Delta}^{F}(E)$ with $F$ of the form 
\begin{equation}\label{eqn F graph}
F=\{(x,\phi_2(x),\dots,\phi_l(x))\, :\, x\in E\},
\end{equation}
 with arbitrary continuous functions $\phi_j:\R^d\to\R^d$. 
 Further comments are in \S\S\ref{subsec relation}. below.
 However, this argument relies on both the form of $F$ and the product nature of the Euclidean metric on $\R^{ld}$,
 and thus does not apply to our most general  result, Thm. \ref{main three}.

\section{The main results}\label{sec main results}

Our main results are the following, in increasing order of generality.

\begin{theorem}\label{main one}
If $E\subseteq \mathbb{R}^d,\, d\geq 2,$ is a Borel set
with $\dim(E) > \frac{2d+1}{3}$, then
 $\Int ( \mathring{\Delta}^{diag} (E)) \neq \varnothing$.  
\end{theorem}

This result is the $3$-point configuration set case of

\begin{theorem}\label{main two}
Let $E\subseteq \mathbb{R}^d$ be compact, $d\ge 2$. 
For $k\ge 3$, define the $k$-point configuration set,
$$\mathring{\Delta}_k^{diag}(E):=\left\{\, \left|\, (x,\dots,x) - \left(y_1,\dots, y_{k-1}\right)\, \right|\,
: \, x, y_1,\dots, y_{k-1}\in E ,\, y_i\ne y_j \right\},$$
$|\, \cdot \, |$ being the Euclidean norm on $\R^{(k-1)d}$.
Then  $\Int ( \mathring{\Delta}_k^{diag} (E)) \neq \varnothing$ if 
\begin{equation}\label{eqn k condition in main two}
\dim(E) > 
\begin{cases}
\frac{2d+1}{3}, & k=3 \\
\frac{(k-1)d}k, & k\ge 4.
\end{cases}
\end{equation}
\end{theorem}

More generally, we have

\begin{theorem}\label{main three}
Suppose $d\ge 2$ and $k\ge 3$.
Let  $\mathcal G$ denote the space of $C^\omega$ Riemannian metrics on $\left(\R^d\right)^{k-1}$.
For any $g\in \mathcal G$,  let $d_g$ be  the induced distance function, which is defined on at least a neighborhood of the diagonal of 
$\left(\R^d\right)^{k-1}$.
Let $g_0$ denote the Euclidean metric.  
Then there is an $N=N_{d,k}\in\mathbb N$ and
a neighborhood $\mathcal U$ of $g_0$ in the $C^N$ topology on $\mathcal G$ such 
that if $g\in\mathcal U$,
and  for a compact $E\subset\R^d$ we define
\beas
\mathring{\Delta}_g^{diag}(E) &:=& \big\{ \, d_g\left(\left(x,\dots,x\right),\left(y_1,\dots,y_{k-1}\right)\right)\, \\
& & \qquad\quad  :\, x,y_1,\dots, y_{k-1} \in E,\, y_i\ne y_j,\, \forall\, i\ne j \big\},
\eeas
then  $\Int ( \mathring{\Delta}^{diag} (E)) \neq \varnothing$ if \eqref{eqn k condition in main two} holds.
\end{theorem}

\subsection{Relations with known results}\label{subsec relation}

As explained in Remark 4 above, a result for the pinned  Mattila-Sj\"olin problem in $\R^d$
automatically yields nonempty interior for $\mathring{\Delta}^F(E)$ whenever $F$ is of the form \eqref{eqn F graph}, 
which includes the $(k-1)$-fold diagonal.
Thus, the pinned distance set threshold of $\dim(E)>(d+2)/2,\, d\ge 3,$ from Peres-Schlag  \cite{PeresSchlag}
produces a better result than Thm. \ref{main one} for $d\ge 4$, and similarly for \cite{IosevichLiu1} for $d\ge 5$.
however, Thm. \ref{main one} is better for $d=3$, and  for $d=2$, where \cite{PeresSchlag} doesn't apply. 
Similarly, \cite{PeresSchlag} yields  for $k\ge 4$  a threshold at least as good as  Thm. \ref{main two} in all  $d\ge 3$.

The very recent paper of Borges,  Iosevich and Ou \cite{BorIosOu1}  gives a lower threshold than our Thm. \ref{main one} in all dimensions.
The authors state that their method extends to the context of Thm. \ref{main two}, but without giving specific thresholds.
On the other hand, it is not clear that the technique of \cite{BorIosOu1} applies in the setting of Thm. \ref{main three},
due to the non-product nature of general Riemannian metrics $g$ on $\left(\R^d\right)^{k-1}$.

It is  reasonable to ask why we are persisting in stating and proving Thms. \ref{main one} and \ref{main two}.
The point is that, rather than proving Thm. \ref{main three} immediately, we will build up to it with a proof of Thm. \ref{main one} 
based on an $L^2 \times L^2 \to L^2$ decay bound for a bilinear spherical averaging operator.
This naturally leads to  the multilinear  operators and estimates yielding Thm. \ref{main two},
which we analyze and prove  in the Fourier integral operator framework  of Greenleaf, Iosevich and Taylor  \cite{GreenleafIosevichTaylor2}.
With minimal additional effort, this then leads to  Thm. \ref{main three} in the case of a product metric;
the inherent stability of the FIO approach under general perturbations then allows it to be proven in full generality.

We now start with the proof of Thm. \ref{main one}.

\section{The Bilinear Spherical Averaging operator}\label{sec bilinear sao}

Let $d\geq 2$. Then for $x\in \mathbb{R}^d$, $t>0$ and for functions $f,g \in \mathcal{S} (\mathbb{R}^d)$ we define the averaging operator:

$$A_r(f,g)(x)= \int_{\mathbb{S}^{2d-1}} f(x-ru)\,g(x-rv)\, d\sigma(u,v)$$
where $\sigma$ is the surface measure on unit sphere $\mathbb{S}^{2d-1}$ in $\mathbb{R}^{2d}$,
$$\mathbb{S}^{2d-1}=\{ (u,v) \in \mathbb{R}^d \times \mathbb{R}^d\, :\, |u|^2+|v|^2=1\, \}$$

Next, we define the (full) maximal version of the bilinear spherical operator

$$\mathcal{M}(f,g)(x)= \sup_{r>0}|A_r(f,g)(x)| $$

Finally, we define its single-scale (localized) bilinear maximal operator, which is

$$\widetilde{\mathcal{M}}(f,g)(x)= \sup_{r\in [1,2]}|A_r(f,g)(x)| $$

\subsection{Known results and goals}

The operators $A_r$ and $\mathcal{M}$ first appeared in the paper of Geba, et al.,  \cite{GGIPES}, where the authors proved some initial $L^p$ improving estimates for these operators. Subsequently, the $L^p$ improving estimates for $\mathcal{M}$ were further developed in the works of Barrionevo-Grafakos-D.He-Honz\'ik-Oliveira (see \cite{BGHHO}), Grafakos-D.He-Honz\'ik (see \cite{GHH}) and Heo-Hong-Yang (see \cite{HHY}). Finally, the full region $L^p$ improving estimates for the operator $\mathcal{M}$ were given in the work of Jeong and Lee (see \cite{JL}) as the result of a clever ``slicing'' argument enabled them to pointwisely dominate the maximal biliner spherical averaging operator by the product of a Hardy-Littlewood maximal operator and a linear spherical averaging operator, both of which have been extensively studied. Furthermore, in the same work the authors explored the $L^p$ improving estimates for the operator $\widetilde{\mathcal{M}}$ getting a large region of exponents, however there is still work left open in this case. Subsequent developments have included sparse domination results \cite{PalssonSovine,BFOPZ} and very recent lacunary maximal operator results \cite{BorgesFoster}.

For our work, we need not only $L^p$ improving estimates but more specifically $L^2\times L^2 \to L^2$ estimate with decay. We already know there operator is bounded from $L^2\times L^2 \to L^2$ but this is not enough. If one considers functions which have compact support on the frequency side, then a decay factor appears. Moreover, we just consider the operator $A_r$ and given the absence of supremum in the definition, one could exploit the full decay of the surface measure on the unit ball in $\mathbb{R}^{2d}$.

We start with the following proposition.

\begin{prop}\label{MainProp}
Let $i, j \in \mathbb{N}$ and let $f,g$ be functions with 

$$supp(\hat{f}) \subset \{\xi \in \mathbb{R}^d \, :\, 2^{i-1}<|\xi| \leq 2^{i+1}\}$$ and 

$$supp(\hat{g}) \subset \{\xi \in \mathbb{R}^d \, :\, 2^{j-1}<|\xi| \leq 2^{j+1}\}$$

then

$$\|A_r(f,g)\|_2 \lesssim_r \big(2^{2i}+2^{2j}\big)^{-\frac{2d-1}{4}} 2^{\min \{i,j\} \frac{d}{2}} \|f\|_2 \|g\|_2$$ 
\end{prop}

\section{Proof of Theorem \ref{main one} } 

In this section we will give a proof of Thm. \ref{main one} using the bilinear spherical averaging operator. 
We start by proving Proposition    \ref{MainProp}.

\subsection{ The decay of the measure \texorpdfstring{$\sigma$}{sig}}

We define $\sigma_r$ to be the surface measure on the sphere of radius $r$ in $\mathbb{R}^{d}\times \mathbb{R}^d$:  

$$\mathbb{S}_r^{2d-1}=\{ (u,v) \in \mathbb{R}^d \times \mathbb{R}^d\, :\, |u|^2+|v|^2=r^2\, \}$$ Then by the classical method of stationary phase we have for $(\xi,\eta) \in \mathbb{R}^d\times \mathbb{R}^d$:

\begin{equation}\label{decay}
\reallywidehat{\sigma_{r}}(\xi,\eta)= r^{2d-1} \reallywidehat{\sigma} (r\xi,r\eta) \text{\quad with \quad } |\reallywidehat{\sigma_{r}}(\xi,\eta)| \lesssim r^\frac{2d-1}{2}\, |(\xi,\eta)|^{-\frac{2d-1}{2}} 
\end{equation}

By a change of variable we get 

$$A_r(f,g)(x)= \frac{1}{r^{2d-1}}\,\int_{\mathbb{S}_r^{2d-1}} f(x-u)\,g(x-v)\, d\sigma_r(u,v)$$

Next, for $\xi\in\mathbb{R}^d$ we use the Fourier inversion formula and Fubini's theorem to write $r^{2d-1}\,\reallywidehat{A_r(f,g)}(\xi)$ as:

\begin{align*}&\int_{\mathbb{R}^d}\, \int_{\mathbb{S}_r^{2d-1}} f(x-u)\,g(x-v)\, d\sigma_r(u,v) \,e^{-2\pi i x \cdot \xi}\, dx = \\
&\int_{\mathbb{R}^d}\, \int_{\mathbb{S}_r^{2d-1}}\, \int_{\mathbb{R}^d}\, \int_{\mathbb{R}^d}  \widehat{f}(\eta)\,\widehat{g}(x')\, e^{2\pi i(x-u)\cdot \eta}\, e^{2\pi i(x-v)\cdot x'} dx'd\eta \, d\sigma_r(u,v) \,e^{-2\pi i x \cdot \xi}\, dx=\\
&\int_{\mathbb{R}^d}\, 
\int_{\mathbb{R}^d}\, \int_{\mathbb{R}^d}  \widehat{f}(\eta)\,\widehat{g}(x')\, e^{2\pi i(x,x,x)\cdot (\eta, -\xi, x')} \,\int_{\mathbb{S}_r^{2d-1}} e^{-2\pi i (u,v)
\cdot(\eta,x')} d\sigma_r(u,v) \, dx' \,d\eta\, dx= \\
&\int_{\mathbb{R}^d}\, \int_{\mathbb{R}^d} \, \widehat{f}(\eta)\,\widehat{g}(x')\, \widehat{{\sigma}_r}(\eta,x') \, \int_{\mathbb{R}^d}\, e^{2\pi i(x,x,x)\cdot (\eta, -\xi, x')} \,dx\, dx'\,d\eta
\end{align*}

and since 

$$\int_{\mathbb{R}^d}\, e^{2\pi i(x,x,x)\cdot (\eta, -\xi, x')} \,dx = \delta(\eta-\xi+x')$$ we get the formula

$$\reallywidehat{A_r(f,g)}(\xi)= \frac{1}{r^{2d-1}}\,\int_{\mathbb{R}^d}\widehat{f}(\eta)\,\widehat{g}(\xi-\eta)\, \widehat{\sigma_{r}}(\eta,\xi-\eta)\, d\eta=\int_{\mathbb{R}^d}\widehat{f}(\eta)\,\widehat{g}(\xi-\eta)\, \widehat{\sigma}(r(\eta,\xi-\eta))\, d\eta$$

\subsection{Proof of Proposition \ref{MainProp}}\label{ProofMain} 
\begin{proof} We can assume without loss of generality that $i\leq j$. Then, we apply Plancherel's theorem to get:

$$\|A_r(f,g)\|^2_2 =
\|\reallywidehat{A_r(f,g)}\|^2_2=
\int_{\mathbb{R}^d} \Big( \int_{\mathbb{R}^d}\widehat{f}(\eta)\,\widehat{g}(\xi-\eta)\, \widehat{\sigma}(r(\eta,\xi-\eta))\, d\eta \Big)^2  d\xi$$

Using the decay of the measure in \eqref{decay} we get 

$$\lesssim_r \big(2^{2i}+2^{2j}\big)^{-\frac{2d-1}{2}} \int_{\mathbb{R}^d} \Big( \int_{\mathbb{R}^d}|\widehat{f}(\eta)|\,|\widehat{g}(\xi-\eta)|\, d\eta \Big)^2 d\xi
$$

since $|\eta|\sim 2^i$, $|\xi-\eta|\sim 2^j$ and so $|(\eta, \xi-\eta)|^2= |\eta|^2+|\xi-\eta|^2 \sim 2^{2i}+2^{2j}$.

Next, let $A^{\xi}_{i,j}:=\{\eta\, : \, |\eta|\sim 2^i\,,\, |\xi-\eta|\sim 2^j \}$. Note that the inner integral is supported on this set and so we can estimate it by using Cauchy-Schwarz inequality:

$$ \int_{A^{\xi}_{i,j}}|\widehat{f}(\eta)|\,|\widehat{g}(\xi-\eta)|\, d\eta  \lesssim 2^{\frac{id}{2}} \Big(\int_{\mathbb{R}^d} \, |\widehat{f}(\eta)|^2 \, |\widehat{g}(\xi-\eta)|^2 d \eta\Big)^{\frac{1}{2}} $$

which gives, after applying Fubini's theorem and a change of variable:

$$\|A_r(f,g)\|^2_2 \lesssim_r \big(2^{2i}+2^{2j}\big)^{-\frac{2d-1}{2}} \, 2^{id} \|f\|_2^2\|g\|_2^2$$

This finishes the proof of Proposition \ref{MainProp}.
\end{proof}

\subsection{Conclusion of proof of Theorem \ref{main one}}\label{subsec conclusion one}
\begin{proof}
Let $E\subset \mathbb{R}^d$ with $\dim(E)>\frac{2d+1}{3}$. Then, fix $s\in (\frac{2d+1}{3}, \dim(E))$. We argue as in the proof of Theorem 4.6 in \cite{Mattila2}. By the Frostman's Lemma \cite[~Thm. 2.8]{Mattila2}, there exists a measure $\mu \in \mathcal{M}(E)$ with $I_s(\mu) <\nolinebreak\infty$. Then we define its \textit{distance measure} $\delta(\mu) \in \mathcal{M}(\mathring{\Delta}^{diag}(E))$ defined for Borel sets $B\subset \mathbb{R}^d$ by

$$\delta(\mu)(B)= \int_{\mathbb{R}^d} \mu\times\mu \big( (y_1,y_2)\, :\, |(x,x)-(y_1,y_2)| \in B \big) d\mu(x)$$

where we seemingly have enlarged the integrand by adding in the case when $y_1=y_2$, but note that $\mu\times\mu((y_1,y_2):y_1=y_2)=0$ which follows from the Frostman condition because the set is of strictly lower dimension than $E\times E$. In other words, $\delta(\mu)$ is the image of $\mu \times \mu$ under the distance map $(x,y_1,y_2) \to |(x,x)-(y_1,y_2)|$, or equivalently for any continuous function $g$ on $\mathbb{R}$:

$$\int_0^{\infty} g(r) d\delta(\mu)(r)= \int_E \,\int_E\,\int_E g\big( |(x,x)-(y_1,y_2)| \big)\, d\mu(x) \, d\mu(y_1)\, d\mu(y_2)$$

Next, for a smooth function $f$ with compact support we have $\delta(f)$ is also a function. To see this we write:

\begin{align*}
\int_{\mathbb{R}} g(r) d\delta(f)(r)&= \\
&=\int_{\mathbb{R}} \,\int_{\mathbb{R}}\,\int_{\mathbb{R}} g\big( |(x,x)-(y_1,y_2)| \big) f(x)\, f(y_1)\,f(y_2)\, dx\, dy_1\, dy_2 \\
&= \int_{\mathbb{R}^{2d}} \int_{\mathbb{R}}  g\big( |(x,x)-y| \big) (f\otimes f)(y)\, f(x)\, dx\,dy \\
&= \int_{\mathbb{R}} \int_{\mathbb{S}^{2d-1}} \int_0^{\infty} g(r) \, f(x-r \omega_1) \, f(x-r\omega_2) \, r^{2d-1}\,dr\, d\sigma(\omega) \, f(x) dx
\end{align*} where in the third equality for $y=(y_1,y_2)\in \mathbb{R}^{d} \times \mathbb{R}^{d} $ we have $(f\otimes f)(y):=f(y_1)\, f(y_2)$. In the third equality we used polar coordinates where $\omega=(\omega_1,\omega_2)\in \mathbb{S}^{2d-1}$. Therefore, after applying Fubini's theorem and using the definition of the Bilinear Averaging Spherical operator $A_r$ we have

$$\int_0^{\infty} g(r) d\delta(f)(r)= \int_0^{\infty} g(r) \,r^{2d-1}\, A_r(f,f)(x) \, f(x) dx dr$$ which implies that $\delta(f)$ is a function with

$$\delta(f)(r)= r^{2d-1}\,\int_{\mathbb{R}^d} A_r(f,f)(x) \, f(x) dx$$ 

Next, we will approximate weakly the measure $\mu$ mentioned above. Namely, let $\psi$ be a smooth compactly supported function in $\mathbb{R}^d$ with $\int \psi =1$. As usually, we define $\psi_{\epsilon}(x)=\epsilon^{-d}\psi(\frac{x}{\epsilon})$ and $\mu_{\epsilon} = \psi_{\epsilon} \ast \mu$. Then $\mu_{\epsilon} \to \mu$, as $\epsilon \to 0$ weakly and so $\delta(\mu_{\epsilon}) \to \delta(\mu)$ weakly. Moreover $\widehat{\mu_{\epsilon}}(x)=\widehat{\psi}(\epsilon\,x) \, \widehat{\mu}(x) \to \widehat{\mu}(x)$ for all $x\in \mathbb{R}^d$. 

As we know, for any $\epsilon>0$, $\mu_{\epsilon}$, is a function. Then, from the formula for $\delta(f)$ we get 

$$ \delta(\mu_{\epsilon})(r)= r^{2d-1}\,\int_{\mathbb{R}^d} \, A_r(\mu_{\epsilon},\mu_{\epsilon})(x) \, \mu_{\epsilon}(x) dx$$ and by the comments above the left side converges weakly to $\delta(\mu)(r)$. We would like to see where the right hand side converges to. Using Parseval's theorem we see 

$$\delta(\mu_{\epsilon})(r)= r^{2d-1}\,\int_{\mathbb{R}^d} \reallywidehat{A_r(\mu_{\epsilon},\mu_{\epsilon})}(\xi) \, \reallywidehat{\mu_{\epsilon}}(\xi) d\xi $$

Next, we have $\lim\limits_{\epsilon\to 0} \reallywidehat{\mu_{\epsilon}}(\xi) =\widehat{\mu}(\xi)$  pointwise, and since $$\reallywidehat{A_r(\mu_{\epsilon},\mu_{\epsilon})}(\xi) = \int_{\mathbb{R}^d} \widehat{\psi}(\epsilon\,\eta \, )\widehat{\mu}(\eta) \,\widehat{\psi}(\epsilon\,(\xi-\eta)) \, \widehat{\mu}(\xi-\eta)\, \widehat{\sigma}(r(\eta,\xi-\eta))\, d\eta$$
we get
$$\lim\limits_{\epsilon\to 0}\reallywidehat{A_r(\mu_{\epsilon},\mu_{\epsilon})}(\xi) = \int_{\mathbb{R}^d} \widehat{\mu}(\eta) \, \widehat{\mu}(\xi-\eta)\, \widehat{\sigma}(r(\eta,\xi-\eta))\, d\eta$$
with passing the limit inside justified by the Dominated Convergence Theorem. Note that

\begin{align*}|\widehat{\mu}(\eta) \, \widehat{\mu}(\xi-\eta)\, \widehat{\sigma}(r(\eta,\xi-\eta))| &\lesssim _r |\widehat{\mu}(\eta) \, \widehat{\mu}(\xi-\eta)||(\eta, \xi-\eta)|^{-\frac{2d-1}{2}} \\
&\lesssim_r |\widehat{\mu}(\eta) \, \widehat{\mu}(\xi-\eta)|\, |\eta|^{-\frac{2d-1}{4}} \, |\xi-\eta|^{-\frac{2d-1}{4}}
\end{align*}
and the last function is integrable (in $\eta$) by utilising the Cauchy-Schwarz inequality, a change of variable and the fact $I_{\frac{1}{2}}(\mu) \leq I_s(\mu)$, since $\mu$ has compact support and $s> \frac{1}{2}$.

Now we write 
$$\delta(\mu_{\epsilon})(r)=r^{2d-1}\,\int_{\mathbb{R}^d} \,B_{\epsilon}(\xi)\, E_{\epsilon}(\xi)\, d \xi$$ 
where 
$$B_{\epsilon}(\xi)= |\xi|^{\frac{d-s}{2}}\int_{\mathbb{R}^d} \widehat{\psi}(\epsilon\,\eta \, )\widehat{\mu}(\eta) \,\widehat{\psi}(\epsilon\,(\xi-\eta)) \, \widehat{\mu}(\xi-\eta)\, \widehat{\sigma}(r(\eta,\xi-\eta))\, d\eta $$
and
$$E_{\epsilon}(\xi) = |\xi|^{\frac{s-d}{2}}\widehat{\psi}(\epsilon\,\xi)\, \widehat{\mu}(\xi)$$

We will dominate each of these functions by $L^2$ integrable functions, independently of $\epsilon$ and so $B_{\epsilon}(\xi)\, E_{\epsilon}(\xi)$ will be 
dominated by an $L^1$ function, independently of $\epsilon$ which will allow us to use the dominated convergence theorem and get the formula

\begin{equation}\label{dform}
\delta(\mu)(r)= r^{2d-1}\int_{\mathbb{R}^d}\Big( \int_{\mathbb{R}^d} \widehat{\mu}(\eta) \, \widehat{\mu}(\xi-\eta)\, \widehat{\sigma}(r\eta,r(\xi-\eta))\, 
d\eta \Big) \widehat{\mu}(\xi) d\xi
\end{equation}

Firstly,
$$|E_{\epsilon}(\xi)| = |\xi|^{\frac{s-d}{2}}|\widehat{\psi}(\epsilon\,\xi)\, \widehat{\mu}(\xi)| \lesssim_{\psi} |\xi|^{\frac{s-d}{2}}\,| \widehat{\mu}(\xi)| $$
and the $L^2$ norm of the right side is exactly equal to $I_s(\mu)$ which is finite. 

Secondly,
  $$|B_{\epsilon}(\xi)| \lesssim_{\psi} |\xi|^{\frac{d-s}{2}}\int_{\mathbb{R}^d} |\widehat{\mu}(\eta)|  \, |\widehat{\mu}(\xi-\eta)|\, |\widehat{\sigma}(r(\eta,\xi-\eta))|\, d\eta$$

Now we will decompose $\hat{\mu}$ on dyadic scales. Consider the Schwartz functions $\eta_0 (\xi)$ supported at $|\xi| \leq \frac{1}{2}$ and $\eta (\xi)$ 
supported in the spherical shell $\frac{1}{2}<|\xi|\leq 2$ such that the quantities $  \eta_0(\xi)$, $  \eta_j(\xi):= \eta \big(2^{-j}\xi\big)$, with $j \geq 1$, form a 
partition of unity.

Then we define $\mu_j(x):= \mu \ast \widecheck{\eta_j}(x)$ and so $ \widehat{\mu}_j(\xi)= \widehat{\mu}(\xi) \eta_j(\xi)$. Thus, $\widehat{\mu}(\xi)=\sum\limits_{j=0}^{\infty} \widehat{\mu_j}(\xi)$ and moreover, 

\begin{align*}
|B_{\epsilon}(\xi)| &\lesssim_{\psi} \sum\limits_{i,j=0}^{\infty} |\xi|^{\frac{d-s}{2}}\int_{\mathbb{R}^d} |\widehat{\mu_i}(\eta)|  \, |\widehat{\mu_j}(\xi-\eta)|\, |\widehat{\sigma}(r(\eta,\xi-\eta))|\, d\eta \\
&\lesssim \sum\limits_{i,j=0}^{\infty} (2^i+2^j)^{\frac{d-s}{2}} \int_{\mathbb{R}^d} |\widehat{\mu_i}(\eta)|  \, |\widehat{\mu_j}(\xi-\eta)|\, |\widehat{\sigma}(r(\eta,\xi-\eta))|\, d\eta 
\end{align*}
since $|\xi| \leq |\eta|+|\xi-\eta|\lesssim 2^i+2^j$ on the supports of $\widehat{\mu_i}, \widehat{\mu_i}$ and  $d>s$. Now the function on the right is independent of $\epsilon$ and $L^2$ integrable as

$$I:=\bigg(\int_{\mathbb{R}^d}\, \Big( \sum\limits_{i,j=0}^{\infty} (2^i+2^j)^{\frac{d-s}{2}} \int_{\mathbb{R}^d} |\widehat{\mu_i}(\eta)|  \, |\widehat{\mu_j}(\xi-\eta)|\, |\widehat{\sigma}(r(\eta,\xi-\eta))|\, d\eta\, \Big)^2 \, d\xi \bigg)^{\frac{1}{2}}  $$

$$\leq \sum\limits_{i,j=0}^{\infty} (2^i+2^j)^{\frac{d-s}{2}} \bigg(\int_{\mathbb{R}^d}\, \Big(\int_{\mathbb{R}^d} |\widehat{\mu_i}(\eta)|  \, |\widehat{\mu_j}(\xi-\eta)|\, |\widehat{\sigma}(r(\eta,\xi-\eta))|\, d\eta\, \Big)^2 \, d\xi \bigg)^{\frac{1}{2}}$$

$$\lesssim_r \sum\limits_{i,j=0}^{\infty} (2^i+2^j)^{\frac{d-s}{2}} \big(2^{2i}+2^{2j}\big)^{-\frac{2d-1}{4}} 2^{\min \{i,j\} \frac{d}{2}} \|\mu_i\|_2 \|\mu_j\|_2$$
by Minkowski's integral inequality and Proposition \ref{MainProp}. Next, we want to evaluate the $L^2$ norms of the functions $\mu_i$. We have by Plancherel's theorem:

\begin{align*}
\|\mu_{i}\|_2^2 =\|\widehat{\mu_{i}}\|_2^2 &\lesssim 2^{i(d-s)} \int_{\mathbb{R}^d} |\xi|^{-d+s} |\widehat{\mu_{i}}(\xi)|^2 d\xi \\
&\lesssim_{\eta} 2^{i(d-s)}\int_{\mathbb{R}^d} |\xi|^{-d+s} |\widehat{\mu}(\xi)|^2 d\xi \\
&= 2^{i(d-s)}I_s(\mu) \\
&\lesssim_{\mu} 2^{i(d-s)}
\end{align*}

With this at hand, we continue estimating $I$ by utilizing the symmetry of the summand

\begin{align*}
I &\lesssim \sum\limits_{i=0}^{\infty} \sum\limits_{j=i}^{\infty} 2^{j\frac{d-s}{2}} \, 2^{-j\frac{2d-1}{2}} \,  2^{i\frac{d}{2}} \, 2^{i\frac{d-s}{2}}\, 2^{j\frac{d-s}{2}} \\
&= \sum\limits_{i=0}^{\infty} 2^{i(d -\frac{s}{2})} \sum\limits_{j=i}^{\infty}  2^{j(\frac{1}{2} -s)}\\
&\lesssim\sum\limits_{i=0}^{\infty} 2^{i(d -\frac{s}{2})} \, 2^{i(\frac{1}{2} -s)}
\end{align*} 
which is finite since $s> \frac{2d+1}{3}$.

Therefore, for a set $E$ with $\dim(E) > \frac{2d+1}{3}$ we have that the function in \eqref{dform} is continuous, as can be seen by the 
Dominated Convergence Theorem. Next, since $supp (\delta(\mu)) \subset {\Delta}^{diag}(supp (\mu)) \subset {\Delta}^{diag}(E)$ it follows that  $\mathring{\Delta}^{diag}(E)$ has non-empty interior.
\end{proof}

\begin{remark}\label{remark 4.1}
The same proof works for an arbitrary number of points. 
Namely, for $k\geq 3$, if for $E\subset\R^d$ compact we define the $k$-point. configuration set
$\mathring{\Delta}_k^{diag}(E)$ as in Thm. \ref{main two},
then 
\begin{center}
$\dim(E) > \frac{(k-1)d+1}{k}$ implies that $\Int ( \mathring{\Delta}_k^{diag} (E)) \neq \varnothing $,
\end{center}
extending what we have just shown for $k=3$.
However, it turns out that by using the Fourier integral operator  (FIO) approach of \cite{GreenleafIosevichTaylor1,GreenleafIosevichTaylor2},
 for $k\ge 4$ one can lower this by $1/k$.
More importantly, the FIO approach does not require the metric to be Euclidean, or a product, or even translation invariant,
leading to Thm. \ref{main three}..
\end{remark}

\section{Thm. \ref{main two} by a Fourier integral operator approach}\label{sec microlocal}

We now prove  Theorem \ref{main two} using multilinear Fourier integral operators (FIO), 
improving on the threshold in Remark \ref{remark 4.1} for $k\ge 4$.
Using the  FIO method developed in  \cite{GreenleafIosevichTaylor1} for 2-point configuration sets
and then extended in \cite{GreenleafIosevichTaylor2} to $k$-point configurations by  optimizing linear FIO estimates over all bipartite partitions  
of the variables,  
this will then set the scene for the proof of Thm. \ref{main three}.

We will give the calculations needed to prove Theorem \ref{main two},
using the general framework and notation of  \cite{GreenleafIosevichTaylor2}, 
which the reader should consult for a full exposition.
In the terminology of   \cite{GGIP12}, the $k$-configuration set $ \mathring{\Delta}_k^{diag} (E)$ is a 
\emph{$\Phi$-configuration set}. 
For convenience, we relabel $(x,y_1,\dots,y_{k-1})$ as $(x^0,x^1,\dots,x^{k-1})$ and define $\Phi:(\R^d)^k\to\R$,

\begin{equation}\label{eqn Phi euclid}
\Phi(x^0,x^1,\dots,x^{k-1})=\frac12\sum_{i=1}^{k-1} \bv x^0-x^i \bv^2,
\end{equation}
so that $\Int ( \mathring{\Delta}^{diag}_k (E)) \neq \varnothing $ iff 
$$\mathring{\Delta}_\Phi(E):=\left\{\Phi\left(x^0,x^1,\dots ,x^{k-1}\right): x^0,\dots, x^{k-1}\in E, x_i\ne x_j,\, \forall i\ne j\right\},$$ 
has nonempty interior.

We start by finding a base point in $\R^{kd}$ about which to work.
 Let $s_0=s_0(d,k)$ be the threshold for $\dim(E)$ in \eqref{eqn k condition in main two} in the statement of Thm. \ref{main two},
and suppose $\dim(E)>s_0$. Pick an $s$ with $s_0<s<\dim(E)$, let $\mu$ be a Frostman measure supported on $E$ and of finite $s$-energy 
(see \cite[~Thm. 8.17]{Mattila}).
We claim that there exist points $x_0^0,\dots,x_0^{k-1}\in E$ and a $\delta>0$ such that $\mu(B(x_0^i,\delta))>0,\, 0\le i\le k-1$, 
$$x_0^i\ne x_0^0,\, \forall\, i>0,\hbox{ and } x_0^i\ne x_0^j,\, \forall\, 1\le i\ne j\le k-1,$$
and then  set 
\begin{equation}\label{eqn tzero}
t_0:=\frac12\sum_{i=1}^{k-1} \left|x_0^0-x_0^i\right|^2>0.
\end{equation}
To see this one can argue as in \cite[~\S\S 4.1]{GreenleafIosevichTaylor2}. The key point is that if we define
\begin{equation}\label{def W}
W:=\left\{(x^0,\dots,x^{k-1})\in\R^{kd}\, :\, x^i\ne x^0,\, \forall\, i>0,\hbox{ and } x^i\ne x^j,\, \forall\, 1\le i\ne j\le k-1\,\right\},
\end{equation}
then $W$ is a Zariski open subset of $\R^{kd}$, whose complement is contained in a union of algebraic varieties of dimensions $\le (k-1)d$
(since each $\{x^i=x^j\}$ is codimension $d$).
Hence $\dim\left(\R^{kd}\setminus W\right)\le (k-1)d<s$, so that $(\mu\times\cdots\times\mu)(\R^{kd}\setminus W)=0$. 
See \cite{GreenleafIosevichTaylor2}, where this type of argument is given for several different $\Phi$-configurations, for more details.

For each $t>0$, the configuration function $\Phi$ induces a surface measure,
$$K_t(\cdot)=\delta\left(\Phi\left(x^0,\dots,x^{k-1}\right) - t \right)\in\mathcal D'\left(\R^{kd}\right),$$
supported on the incidence relation
\begin{equation}\label{eqn Z euclid}
Z_t:=\left\{\, \left(x^0,\dots,x^{k-1}\right)\in\R^{kd}\, :\, \Phi(x^0,\dots,x^{k-1})=t\,\right\}.
\end{equation}
$K_t$ is a Fourier integral distribution on $\R^{kd}$,
$$K_t(\cdot)=\int_\R e^{i\tau\left(\Phi(\cdot)-t\right)} 1(\tau)\, d\tau\in I^{0+\frac12-\frac{kd}4}\left(N^*Z_t\right),$$
where $N^*Z_t\subset T^*\left(\R^{kd}\right)\setminus 0$ is the conormal bundle of $Z_t$. For convenience, we will write $N^*Z_t$ 
with each pair of spatial and cotangent variables, $(x^i,\, \xi^i)$, grouped together. 
Thus, 
\begin{eqnarray}\label{eqn nzt euclid}
N^*Z_t &=& \big\{ \big(x^0,\tau\sum_{i=1}^{k-1}(x^0-x^i); x^1,-\tau(x^0-x^1); \dots; x^{k-1},-\tau(x^0-x^{k-1})\big) \nonumber\\
& & \qquad :\, \left(x^0,\dots,x^{k-1}\right)\in Z_t,\, \tau\ne 0 \, \big\}.
\end{eqnarray}
To make this more explicit, we parametrize an open subset of $Z_t$ by letting $x^0$ range freely over $\R^d$, and then write $x^i=x^0+y^i,\, 1\le i\le k-1$.
Writing  $\vec{y}=(y^1,\dots,y^{k-2})\in\R^{(k-2)d}$, set 
$r(\vec{y},t)=\left( 2t -  \sum_{i=1}^{k-2} |y^i|^2\right)^\frac12$ and let
\beas
\mathring{U}_t&:=& \big\{ \, (\vec{y},y^{k-1})\in \R^{(k-1)d}:\, 
|y^i|>0,\,\forall\, 1\le i\le k-1; \\
& & \qquad\quad \sum_{i=1}^{k-2} |y^i|^2<2t;  y^{k-1} =r(y^1,\dots,y^{k-2},t)\omega,\, \omega\in\sph;\\
& & \qquad\quad \hbox{ and } y^{i}\ne y^j,\, \forall\, 1\le i\ne j\le k-1 \big\},
\eeas
which is an open subset  of $\R^{(k-1)d}$. 
Since all of the $x^i-x^0=y^i$ are 
distinct, it follows that  $x^i\ne x^j,\, \forall\, 0\le i\ne j\le k-1$. Thus, 
\beas
Z_t\supset \mathring{Z}_t&:=& \big\{\left(x^0,x^0+y^1,\dots,x^0+y^{k-2},x^0+y^{k-1}\right)\, \\
& & \qquad\qquad  :\, x^0\in\R^d,\,(\vec{y},y^{k-2})\in\mathring{U}_t,\, \big\},
\eeas
allowing us to parametrize the open subset $N^*\mathring{Z}_t\subset N^*Z_t$ as
\begin{eqnarray}\label{eqn nzt}
N^*\mathring{Z}_t&=&\big\{\, \big(x^0,-\tau\left(\sum_{i=1}^{k-2} y^i + r(\vec{y},t)\omega\right); x^0+y^1,\tau y^1;\dots;\nonumber \\
& & \qquad\quad x^0+y^{k-2},\tau y^{k-2}; x^0+r(\vec{y},t)\omega,\tau r(\vec{y},t)\omega\big) \\
& & \qquad\quad\, :\,  x^0\in\R^d,\, (\vec{y},y^{k-1})\in \mathring{U}_t,\, \omega\in\sph,\, \tau\ne 0\, \big\}.\nonumber
\end{eqnarray}

Note that $(x_0^0,x_0^1,\dots,x_0^{k-1})\in \mathring{Z}_{t_0}$, with $t=t_0$ as in \eqref{eqn tzero}.
Multiplying $K_t$ by a smooth cutoff function in order to localize to where $x^i\ne x^j,\, \forall\, 0\le i\ne j\le k-1$,
yields an element of $I^{\frac12-\frac{kd}4}\left(N^*\mathring{Z}_{t}\right)$, which for simplicity we still denote by $K_{t}$.

With all of this in place, we commence the proof of Thm. \ref{main two}  by treating the case $k=3$,  showing  how the FIO approach 
reproves Thm. \ref{main one}. Since the number of points is odd,
we need to consider partitions of the form $\sigma=\left(\, \sigma_L\, |\, \sigma_R\, \right)= (\, i\, |\, j\, k\, )$,
with $i,j,k\in\{0,1,2\}$ distinct, and in fact we focus on $\sigma_0:=(\,0\, |\, 12\,)$. 
This corresponds treating $K_{t_0}$ as the Schwartz kernel
of a linear FIO, $T_{t_0}^{\sigma_0}$,  taking functions of $x^1,x^2$ to functions of $x^0$. 
From \eqref{eqn nzt}, the canonical relation $C_t^{\sigma_0}$ for general $T_t^{\sigma_0}$ simplifies to
\begin{eqnarray}\label{eqn ct k is three}
C^{\sigma_0}_t&=&\big\{\, \big(x^0,-\tau\left( y^1+ r(y^1,t)\omega\right); x^0+y^1,\tau y^1; x^0+r(y^1,t)\omega,\tau r(y^1,t)\omega\big) \nonumber \\
& & \qquad\qquad :\,  x^0\in\R^d,\, (y^1,y^2)\in \mathring{U}_t\subset\R^{2d},\, \omega\in\sph,\, \tau\ne 0\, \big\},
\end{eqnarray}
where $r(y^1,t)=\left(2t-|y^1|^2\right)^\frac12$. 
$C^{\sigma_0}_t$ has dimension $3d$ and 
$C^{\sigma_0}_t\subset \left(T^*\R^d\setminus 0\right) \times \left(T^*\R^{2d}\setminus 0\right)$,
so that  $T_t^\sigma:\mathcal E'(\R^{2d})\to \mathcal D'(\R^d)$. We claim that $C_t^{\sigma_0}$ is  \emph{nondegenerate},
in the sense that the projections \linebreak$\pi_L:C_t^{\sigma_0}\to T^*\R^d$ and $\pi_R:C_t^{\sigma_0}\to T^*\R^{2d}$ have maximal rank, i.e., 
are a submersion and an immersion, resp.
By a general principle  concerning canonical relations, the two projections have the same corank at each point, 
so we only need to verify this for $\pi_L:C^\sigma_t\to T^*\R^{d}$,
which follows from $|D(x^0,\xi^0)/D(x^0,\omega,\tau)|\ne 0$. By \cite[~Thm. 5.2(ii)]{GreenleafIosevichTaylor2} (with $p=1$ since $\Phi$ is $\R^1$-valued),
$$\hbox{ if } \dim(E)>\frac13[\max({d,2d}+1)]=\frac{2d+1}3,\hbox{ then } \Int\left(\mathring{\Delta}_\Phi(E)\right) \ne \varnothing,$$
as desired. 

One can check that the other choices  of $\sigma$, up to irrelevant permutations, namely $(1|02)$ and $(2|01)$, could also have been used, but do no better. For $k\ge 4$ below, we will again exhibit one partition that implies the claimed threshold.
However, when $k\ge 4$, the geometry of the $C_t^\sigma$ is less favorable than for $k=3$: the only way to partition the variables to make $C_t^\sigma$ nondegenerate is to make the total spatial dimension $d_L$ of the variables on the left much less than the dimension $d_R$ on the right,
and then that incurs a penalty by raising the effective order of the associated linear FIO, $T_t^\sigma$. 
(See the discussion in \cite[Sec. 5]{GreenleafIosevichTaylor2}.) On the other hand, for $d_L$ as close to $d_R$ as possible, 
we will see that the projections drop rank,   resulting in $T_t^\sigma$ losing derivatives on $L^2$-based Sobolev spaces. 
This forces us to use \cite[~Thm. 5.2(i)]{GreenleafIosevichTaylor2} in place of part (ii). 

To start, suppose   that $k\ge 4$ is even.
Partition the variables $x^0,x^1,\dots,x^{k-1}$ into two groups of equal cardinality $k/2$ on the left and right, resp., picking
$$\sigma=\big(\sigma_L\, |\, \sigma_R\,\big)=\Big(\, 0\, 1 \, \cdots \frac{k-2}2\,\, \big|\, \frac{k}2 \,\, \cdots\,  k-1\, \Big).$$
Using \eqref{eqn nzt}, one sees that this choice has the following properties.

(i) The total spatial dimensions on the left and right groups  are $d_L=d_R=kd/2$.

(ii)  Using $\sigma$ to rearrange $N^*\mathring{Z}_t$ into $C_t^\sigma\subset T^*\R^{kd/2} \times T^*\R^{kd/2}$, 
\linebreak $C_t^\sigma$ avoids the zero sections on both sides, so that $K_t$ is the Schwartz kernel of a linear Fourier integral operator,
$T_t^\sigma\in I^{\frac12-\frac{kd}4}(\R^{\frac{kd}2},\R^{\frac{kd}2};\, C_t^\sigma)$.

(iii) $C_t^\sigma$ has the property that the projections to the left and right have rank at least $(k+2)d/2+1$.
As remarked above, we only need to verify this for $\pi_L:C^\sigma_t\to T^*\R^{kd/2}$. 
From \eqref{eqn nzt},
one calculates that  $D\pi_L$ restricted to
$$span\left\{\, T_{x^0}\R^d,\, \left\{T_{y^i}\R^d\right\}_{i=1}^{(k-2)/2},\, y^{\frac{k}2}\cdot\partial_{y^{\frac{k}2}},\, T_\omega\sph,\,  \partial_\tau\right\}$$
is injective.  
This uses the fact that the radial derivative of $r(\vec{y},t)$ with with respect to $y^{\frac{k}2}$ is nonzero.
(We could have used $y^i\cdot\partial_{y^i}$ for any of the variables $y^i,\, k/2\le i\le k-2$.)

(iv) Since $\rank(D\pi_L)\ge (k+2)d/2+1$ at each point of $C_t^\sigma$, it follows that $\corank(D\pi_L)=kd-\rank(D\pi_L)\le (k-2)d/2 -1$.
By H\"ormander's estimates for FIO, $T_t^\sigma$ loses at most $\beta^\sigma=(k-2)d/4 -1/2$ derivatives on $L^2$-based Sobolev spaces,
and this is locally uniform in $t$. (To use these estimates, one also needs that the spatial projections
from $C_t^{\sigma}$ onto the left and right variables are submersions, which is easily verified.)

(v) In the notation of \cite[~Thm. 5.2(i)]{GreenleafIosevichTaylor2}, $d_L=d_R=kd/2,\, p=1$, $2\beta^\sigma \le (k-2)d/2 -1$ and 
all of the $E_i=E$;
hence, by that result, $\Int\left(\mathring{\Delta}_k^{diag}(E)\right)\ne\varnothing$ if

\beas
\dim(E)&>& \frac1k\left( \max(d_L,d_R)+p+2\beta^\sigma\right)\\
&=&\frac1k\left(\frac{kd}2+1+\frac{(k-2)d}2-1\right)=\frac{(k-1)d}k,
\eeas
finishing the proof of Thm. \ref{main two} for $k\ge 4$ and even.
\medskip

Finally, for $k\ge 5$ and odd, parity prevents the existence of an equidimensional partition, as it did for $k=3$. Choosing
$$\sigma=\big(\sigma_L\, |\, \sigma_R\,\big)=\Big(\, 0\, 1 \, \cdots \frac{k-3}2\,\, \big|\, \frac{k-2}2 \,\, \cdots k-1\, \Big),$$
one has $d_L=(k-1)d/2 < d_R=(k+1)d/2$, and the resulting canonical relation $C_t^\sigma\subset T^*\R^{(k-1)d/2} \times \R^{(k+1)d/2}$ avoids the zero section on each side.
Calculations almost identical to the case of $k$ even  show that $\rank(D\pi_L)\ge (k+1)d/2+1$ at each point of $C_t^\sigma$,
so that $\corank(D\pi_L)=(k-1)d-\rank(D\pi_L)
\le (k-3)d/2-1$ and thus \cite[~Thm. 5.2(i)]{GreenleafIosevichTaylor2}, 
with $\max(d_L,d_R)=(k+1)d/2,\, p=1, \,  2\beta^\sigma = (k-3)d/2-1$, implies that
 $\Int\left(\mathring{\Delta}_k^{diag}(E)\right)\ne\varnothing$ if
$$\dim(E)>\frac1k\left(\frac{(k+1)d}2+1+\frac{(k-3)d}2-1\right)=\frac{(k-1)d}k,$$
finishing the proof of Thm. \ref{main two}.

\section{Riemannian setting: proof of Thm. \ref{main three}}\label{sec riemann}

H\"ormander's estimates for FIO on $\R^d$ or a $d$-dimensional smooth manifold are stable with respect to perturbations (in $C^\omega$) 
which are small in the $C^N$ topology
on the canonical relations and amplitude, for some $N=N_d$. (This might not be stated explicitly in the literature, but is
a folk theorem, being clear from the proofs.)
Due to this stability,  Thm. \ref{main three}  follows
almost immediately from the proof of Thm. \ref{main two} above.
Perturbing the Euclidean metric $g_0$ on $\R^{(k-1)d}$ in the $C^{N+3}$ topology (for $N=N_{(k-1)d}$)
results in a $C^{N+1}$ perturbation of the geodesic flow,
and hence a $C^{N+1}$ perturbation of  the distance function.
Thus,  the configuration function,
$$\Phi_g(x^0,\dots,x^{k-1})=\frac12d_g\left(\, (x^0,\dots,x^0)\, ,\, (x^1,\dots,x^{k-1})\, \right)^2$$
is a $C^{N+1}$ perturbation of 
$$\Phi_{g_0}(x^0,\dots,x^{k-1})=\frac12\left|\, (x^0,\dots,x^0) - (x^1,\dots,x^{k-1})\, \right|^2,$$
which was the starting point \eqref{eqn Phi euclid} for the analysis in the previous section.
The existence of a base point $(x_0^0,x_0^1,\dots, x_0^{k-1})\in\R^{kd}$ around which to run the whole argument
follows as before, since the analogue of $W$ in the Riemannian case of  \eqref{def W} from the Euclidean case
is again an analytic variety of codimension $\ge d$ and thus has measure zero w.r.t. $\mu\times\cdots\times\mu$.
Forming $\mathring{Z}_t^g$ as above, it is a $C^{N+1}$ perturbation of  $Z_t^{g_0}$
and hence the conormal bundle $N^*Z_t^g$ is a $C^N$ perturbation of \eqref{eqn nzt euclid}.
As a result, for the same choices of partitions $\sigma$ as in the Euclidean case,
the canonical relations $C_t^\sigma$ in the Riemannian case are $C^N$ perturbations of the $C_t^\sigma$ analyzed in the previous section.
Since $C^N$ perturbations of submersions are submersions, this means the same $L^2$-Sobolev estimates hold, yielding nonempty interior of 
$\mathring{\Delta}_g^{diag}(E) $ for the same lower bounds on $\dim(E)$  \eqref{eqn k condition in main two} as in the Euclidean case.

\section{Acknowledgments}
JG was supported in part by NSF grant DMS-1907435;
AG  by NSF grant DMS-2204943; 
EP  by Simons Foundation Grant \#360560; and 
GP  by an AMS-Simons Travel Grant. EP would like to thank the Vietnam Institute for Advanced Study in Mathematics (VIASM) for the hospitality and for the excellent working condition, where part of this work was done.

\end{document}